\documentclass[12pt]{amsart}

\usepackage{amsmath,amsthm,amsfonts,latexsym,amssymb,enumerate,color}
\usepackage{graphicx}
\usepackage{enumitem}

\newcommand\CC{{\mathbb C}}
\newcommand\GA{\mathcal A}

\newcommand\DD{{\mathbb D}}

\newcommand\NN{{\mathbb N}}

\newcommand\TT{{\mathbb T}}

\newcommand\Id{{\rm Id}}
\def\Im{\mathop{\rm  Im}\nolimits}

\def\beq{\begin{equation}}
\def\eeq{\end{equation}}
\newtheorem{thm}{Theorem}[section]
\newtheorem{prop}[thm]{Proposition}

\newtheorem{cor}[thm]{Corollary}
\newtheorem{rem}[thm]{Remark}

\newcommand\beginpf{\noindent {\bf Proof:} \quad}

\newcommand\re{\mathop{\rm Re}\nolimits}
\newcommand\im{\mathop{\rm Im}\nolimits}
\def\beginpf{\begin{proof}}
\def\endpf{\end{proof}}

\renewcommand\phi{\varphi}

\newcommand\1{{\bf 1}}
\newcommand\DW{\mathop{\rm DW}}
\newcommand\fix{\mathop{\rm fix}}
\newcommand\Hol{\mathop{\rm Hol}}
\newcommand\LL{\mathcal L}
\newcommand\T{T_{w,\phi}}

\begin{document}
\title[Weighted composition operators]{Weighted composition operators: isometries and asymptotic behaviour
}

\author{I. Chalendar}
\address{Isabelle Chalendar, Universit\'e Paris Est, 
LAMA, (UMR 8050), UPEM, UPEC, CNRS,
F-77454, Marne-la-Vall\'ee (France)}
\email{isabelle.chalendar@u-pem.fr}

\author{J.R. Partington}
\address{Jonathan R. Partington, School of Mathematics, University of Leeds, Leeds LS2 9JT, UK}
 \email{J.R.Partington@leeds.ac.uk}
 
 %\noindent\textsc{Mathematics Subject Classification} (2000):
 %Primary:  47A15, 47D03
 %Secondary:  30H10, 31C25
 
\baselineskip18pt

\bibliographystyle{plain}

\begin{abstract}
This paper studies the behaviour of iterates of weigh\-ted composition operators acting on spaces of 
analytic functions, with particular emphasis on the Hardy space $H^2$. Questions relating to 
uniform, strong and weak convergence are resolved in many cases.
Connected to this
is the question when a weighted composition operators is an isometry, and new
results are given in the case of the Hardy and Bergman spaces.
\end{abstract}

 \maketitle

\section{Introduction and notation}

Let $H^2$ denote the Hardy space of the disc, and for $w \in H^2$ and $\phi: \DD \to \DD$ holomorphic,
we consider the ({\em a priori\/} densely defined) weighted composition operator $T_{w,\phi}$ on $H^2$
given by
\[
(T_{w,\phi}f)(z)= w(z) f( \phi(z)) \qquad (f \in H^2, \quad z \in \DD).
\]

In general it is hard to determine when $\T$ is bounded  but, 
by  Littlewood's subordination theorem, $w \in H^\infty$ is always a sufficient condition.
The extreme case that   $\T$ is bounded  if and only if $w \in H^\infty$ occurs precisely when $\phi$ is a finite Blaschke product  \cite{CH01,matache08}.
At the other extreme, $w \in H^2$ is necessary and sufficient for boundedness if and only if $\|\phi\|_\infty<1$ \cite{GKP10}.

Continuing   some of the work of \cite{ACKS}, which considered composition operators $C_\phi$
defined by $(C_\phi f)(z)=f(\phi(z))$ for $z \in \CC$,
we shall be considering the following properties that powers of $\T$ may possess:
\begin{itemize}
\item uniform convergence, convergence in the operator norm;
\item strong convergence, the norm convergence of $\T^n f$ for all $f \in H^2$;
\item weak convergence, the weak convergence of $\T^n f$ for all $f \in H^2$ (implying 
local uniform convergence,  the convergence of $\T^n f$ in the Fr\'echet space $\Hol(\DD)$
of holomorphic functions on $\DD$).
% \item ergodicity, the existence of $\lim_{n \to \infty} \frac1n \sum_{k=0}^{n-1} \T^k f$ for all $f \in H^2$.
\end{itemize} 

Each of the properties above implies those below it.

Using the notation $\phi_k$ for the $k$th iterate of $\phi$, 
we have the following easily-derived formula for $\T^n$ for $n=2,3,\ldots$:
\beq\label{eq:tnf}
(\T^n f)(z)=w(z)w(\phi(z))\ldots w(\phi_{n-1}(z))f(\phi_n(z)).
\eeq

Recall that every holomorphic mapping $\phi: \DD \to \DD$ that is not an automorphism
with a fixed point in the disc possesses a {\em Denjoy--Wolff point}, $\DW(\phi)  \in \overline \DD$
such that the iterates of $\phi$ tend to $\DW(\phi)$ uniformly on compact subsets of $\DD$. 
(See, for example \cite[Chap.~2]{CM95}.)

As regards $H^2$, some of the results of \cite{ACKS} may be summarized as follows:
\begin{itemize}
\item Theorem 3.4. $(C^n_\phi)$ converges in operator norm if and only if the essential
spectral radius satisfies $r_e(C_\phi)<1$.
\item Theorem 4.3. If $\DW(\phi) \in \TT$, then $\sup_n \|C^n_\phi\|=\infty$ and so $(C^n_\phi)$
cannot even converge weakly.
\item Theorem 4.8, Theorem 4.10. If $\DW(\phi)=\alpha \in\DD$ and $\phi$ is inner, then $(C^n_\phi f)$ converges
weakly to $P_\alpha f:= f(\alpha) {\bf 1}$ for all $f \in H^2$, where $\bf 1$ is the constant function;
however, it does not converge
strongly.
\item Corollary 4.16. If $\DW(\phi)=\alpha \in \DD$ and $\phi$ is not inner, then $(C^n_\phi)$ converges in operator norm to $P_\alpha$.
\end{itemize}

For mappings $\phi$ with a fixed point $a\in \DD$, we may study the iterates
of $T_{w,\phi}$, reducing to the case where the fixed point is $0$ by using the involution 
\beq\label{eq:psi}
\psi_a: z  \mapsto \frac{a-z}{1-\overline a z}.
\eeq
This provides an isomorphism $C_{\psi_a}$ of $H^2$ and satisfies
\beq\label{eq:conj}
(C_{\psi_a} \T C_{\psi_a} f)(z)= w(\psi_a(z)) f(\psi_a \circ \phi \circ \psi_a (z)),
\eeq
where now $\psi_a \circ \phi \circ \psi_a$ has its fixed point at the origin.\\

In the next few sections we analyse convergence for various classes of $\phi$. 
As will be seen, some of our results
apply to other spaces of functions, although our main focus is on $H^2$.\\

We write $\sigma(T)$ and $\sigma_p(T)$ for the spectrum and point
spectrum (eigenvalues), respectively, of an operator $T \in \LL(X)$, the 
algebra of bounded operators on a Banach space $X$.
Also $r_e(T)$ denotes the essential spectral radius of $T$.

\section{Convergence criteria for weighted composition operators}
\label{sec:2}

There is a general convergence criterion for convergence of powers of an operator $T \in \LL(X)$,
where $X$ is a Banach space, which may be found in \cite[A-III, Sec. 3.7]{crowd}; it was
applied to composition operators in \cite{ACKS}, and we shall also make use of it in this section.

An operator $T \in \LL(X)$ is power-bounded if $\sup_{n\geq 1}\|T^n\|<\infty$. 
 
\begin{prop}\label{prop:2.1}
Let $T \in \LL(X)$  be power-bounded.  Then the following assertions are equivalent:\\
(i) $P:=\lim T^n$ exists in $\LL(X)$ and is of finite rank.\\
(ii) (a) $r_e(T)<1$;\\
(b) $\sigma_p(T) \cap \TT \subseteq \{1\}$;\\
(c) if $1$ is in $\sigma(T)$ then it is a pole of the resolvent of order 1.\\
In that case, $P$ is the residue at $1$.
\end{prop}

In this section we take a Banach space $X$ of holomorphic functions on $\DD$
containing the constant function $\1$, such that there is
a continuous injection $X \to \Hol(\DD)$, and hence point evaluations in $\DD$ are
continuous functionals on $X$. There are numerous examples, including
the Hardy spaces $H^p$ and weighted Bergman spaces $A^p_\alpha$ for $p \ge 1$ and $\alpha> -1$.
We also suppose that $\T$ is a power-bounded
weighted composition operator on $X$, excluding the case when $\phi$ is an elliptic automorphism.

\begin{thm}\label{thm:2.2}
With $X$ and $\T$ as above, 
suppose also that $\alpha:=\DW(\phi) \in \DD$. Then
the sequence $(\T^n)$ converges weakly as $n \to \infty$
if and only if
(i)~$|w(\alpha)|<1$ or $w(\alpha)=1$,
and (ii)~$\sup_n \|\T^n\|< \infty$.
\end{thm}
\beginpf
(i) Since $\T^n(\1)(z)=w(z)\ldots w(\phi_{n-1}(z))$, it follows that  \[\T^n(\1)(\alpha)=w(\alpha)^n.\]
 So we have $|w(\alpha)|<1$ or $w(\alpha)=1$, as well as the power-boundedness condition, if $(\T^n)$ converges weakly.
\\
(ii) Conversely, if $|w(\alpha)|<1$ or $w(\alpha)=1$, then we may first suppose without loss of generality that
$\alpha=0$, by applying the transformation given in \eqref{eq:conj}.
We also have by Schwarz's lemma that under these circumstances $|\phi(z)| \le \delta |z|$ on compact discs
$r\DD$ for $0 < r < 1$, for some
fixed $\delta \in (0,1)$, and so $|\phi_k(z)| \le \delta^k |z|$ for $k=1,2,\ldots$.
Thus $f(\phi_n(z)) \to 0$ for $f \in X$. 

If $|w(0)|<1$, then clearly $w(z)w(\phi(z))\ldots w(\phi_{n-1}(z)) \to 0$ uniformly on any compact $K \subset \DD$.

If $w(0)=1$, then $|w(\phi_{k}(z))-1| \le C \epsilon^{k}$ uniformly on $K$ for some fixed $C>0$ and $\epsilon \in (0,1)$, depending on 
the modulus of $w'$. This ensures the
convergence of the infinite product
$w(z)w(\phi(z))\ldots w(\phi_{n-1}(z))$ to a function $\tilde w(z))$. Since the 
$\T^n$ are uniformly bounded in norm, this implies the weak convergence of $(\T^n)$.
\endpf

We may now apply Proposition \ref{prop:2.1} to obtain the following result.

\begin{thm}\label{thm:2.3}
Under the hypotheses of Theorem \ref{thm:2.2}, suppose that 
either $|w(\alpha)|<1$ or $w(\alpha)=1$,
%$w(\alpha) \in \{0,1\}$ 
and that $\sup_n \|\T^n\|<\infty$.
Then
$(\T^n)$ converges uniformly if and only if $r_e(\T)<1$.
\end{thm}
\beginpf
Once more we may assume without loss of generality  that $\alpha=0$.
Given that $(\T^n)$ converges weakly, we also have the 
following condition (D):
\beq
\label{eq:d}
X = \fix(\T) \oplus \overline{\im(I-\T)}
\eeq
where $\fix(\T)=\{f \in X: f(z)=w(z)f(\phi(z))\}$
as in \cite{ACKS}.

%Now if $w(0)=0$ and $\T f=\lambda f$, we have
%\[
%\lambda^n f(z)= \T^n f(z) = w(z) w(\phi(z))\ldots w(\phi_{n-1}(z)) f(\phi_n(z)),
%\]
%and therefore if $\lambda \ne 0$, $f$ has a zero of order at least $n$ at $0$. So $f$ is identically zero.
%Thus $\fix(\T)=\{0\}$ and $\T$ has no point spectrum on the circle.

Now if $|w(0)|<1$ and $\T f=\lambda f$, we have
\[
\lambda^n f(z)= \T^n f(z) = w(z) w(\phi(z))\ldots w(\phi_{n-1}(z)) f(\phi_n(z)),
\]
and therefore if $|\lambda| \ge 1$ and $|z|<1$, we have 
\[
f(z) = \lim_{n \to \infty} \lambda^{-n} w(z) w(\phi(z))\ldots w(\phi_{n-1}(z)) f(\phi_n(z)) = 0.
\]
Thus $\T$ has no point spectrum on the circle. 

Suppose now that $w(0)=1$. If $f \in \fix(\T)$ we  have
\[
f(z)=\lim_n w(z) w(\phi(z)) \ldots w(\phi_{n-1}(z)) f(\phi_n(z)) ,
\]
and so $\fix(\T)$ is one-dimensional and spanned by 
\[
\widetilde w(z):=\lim_n w(z)w(\phi(z))\ldots w(\phi_{n-1}(z)).
\]

Next, if $\T f=\lambda f$ with $\lambda \in \TT \setminus \{1\}$; then we have
$  w(z)  f(\phi(z))=\lambda f(z)$, and taking $z=0$ we have that $f(0)=\lambda f(0)$ so $f(0)=0$.

We may then write $f(z)=z f_1(z)$ for some $f_1 \in \Hol(\DD)$, and likewise $\phi(z)=z \phi_1(z)$, to
obtain
$w(z) z\phi_1(z) f_1(\phi(z)) = \lambda z f_1(z)$, so that $\phi_1(0)f_1(0)=\lambda f_1(0)$. But $|\phi_1(0)|<1$,
since $\phi_1$ is not an automorphism of $\DD$. So $f_1(0)=0$ and we may write $f_1(z)=zf_2(z)$
for some $f_2 \in \Hol(\DD)$.

We now have
$w(z) \phi(z)^2 f_2(\phi(z)) = \lambda z^2 f_2(z)$, giving $\phi_1(0)^2 f_2(0) = \lambda f_2(0)$
and now $f_2(0)=0$.

Continuing in this way we conclude that $f$ is identically zero, and so $\sigma_p(T) \cap \TT \subseteq \{1\}$
in both the  cases $w(0)=0$ and $w(0)=1$; moreover,  if $r_e(\T)<1$ then $\Id-\T$ is Fredholm, and
so $\Im(\Id-\T)$ is closed.

Moreover, since $\ker (\Id-\T)$ is either $0$ or $1$-dimensional, we see that $1$ is a pole of the
resolvent of order at most 1.

We now have all the conditions of Proposition \ref{prop:2.1}, and may conclude that $\|\T^n-P\| \to 0$,
where the projection $P$ is given by 
% $P=0$ if $w(0)=0$ and 
$P=0$ if $|w(0)|<1$ and
$Pf(z)=\widetilde w(z) f(0)$ if $w(0)=1$.

\endpf

Some applications of these results will be given in the next section.

\section{Automorphisms}
\label{sec:3}

As mentioned above, in the case that $\phi$ is a finite Blaschke product, a necessary and sufficient condition
for boundedness of $\T$ on $H^2$ is that $w \in H^\infty$. We consider some special cases of this situation
(we return to looking at operators on more general spaces in Subsection~\ref{sec:6}).

% \subsection{Automorphisms}

We begin with the case that $\phi$ is an elliptic automorphism of finte order.

Suppose that $\phi(z)= \lambda z$, where $\lambda$ is a primitive $k$th root of unity. In
this situation an important role is played by the function $v \in H^\infty$ defined by
\begin{equation}
\label{eq:defv}
v(z)=w(z) w(\lambda z) \cdots w(\lambda^{k-1}z) \qquad (z \in \DD).
\end{equation}]
We have that 
\[
(\T^{mk+j} f)(z)= v(z)^m w(z)w(\phi(z))\ldots w(\phi_{j-1}(z))f(\lambda^j z).
\]
for $m \ge 0$ and $0 \le j < m$.

\begin{thm}
Suppose that $\phi(z)= \lambda z$, where $\lambda$ is a primitive $k$th root of unity, and let $v$ 
be as in \eqref{eq:defv}.
\begin{itemize}	
\item[(i)] If $\|v\|_\infty < 1$, then $\|T^n_{w,\phi}\| \to 0$ as $n \to \infty$.
\item[(ii)] If $\|v\|_\infty > 1$, then $T^n_{w,\phi}$ does not converge weakly.
\item[(iii)] If $\|v\|_\infty=1$ and $v$ is non-constant, then
 \begin{itemize}
\item[(a)] if $|v|<1$ a.e., it follows that $T^n_{w,\phi} \to 0$ strongly, although not uniformly;
\item[(b)] if $|v|=1$ on a set of positive measure, then $\|v^n\|_2 \not\to 0$ and so 
there is weak convergence of  $T^n_{w,\phi}$ to $0$, but not strong convergence.
\end{itemize}
\end{itemize}
\end{thm}
\beginpf
The proof is made simpler by observing that, apart from the factor $v^m$, there
are only finitely many possibilities for the other factor in $T^n_{w,\phi}f(z)$, and
we may consider subsequences determined by $n=mk+j$ for a fixed  $j$.

Then (i) is clear, and (ii) follows on noting that there is a $z_0 \in \DD$ with $|v(z_0)| > 1$
and examining the values at $(T^n_{w,\phi}f)(z_0)$.

For (iii) we note that in case (a) $v^n \to 0$ a.e. on $\TT$, and the strong convergence follows,
by dominated convergence,
although $\|v^n\|_\infty \not\to 0$. Finally in case (b), weak convergence holds since
$|v(z)|<1$ for all points $z \in \DD$, but strong convergence to $0$ is clearly not satisfied.

\endpf

%Less can be said when $\phi$ is an elliptic automorphism of in finite order.
%Here we take $\phi(z)= \lambda z$, where $|\lambda|=1$ and $\lambda$ is not a root of unity.
%
%It is convenient to write $w=w_i w_o$ (an inner--outer factorization), and consider
%first the behaviour of $w_o(z)w_o(\phi(z))\ldots w_o(\phi_{n-1}(z))$. 
%
%Note that we can write $w_o=\exp u$, where $u \in H^1$. Indeed, to within an additive constant in $i\RR$,
%\[
%u(z)=\frac{1}{2\pi} \int_0^{2\pi} \frac{e^{it}+z}{e^{it}-z}\log |w_o(e^{it})| \, dt.
%\]
%Here we recall the pointwise ergodic theorem   for $f \in L^1(\TT)$.
%Noting that the transformation $z \to \lambda z$ is ergodic (it fixes a.e.\ no subsets of the
%circle of normalized measure strictly between 0 and 1), we have
%\[
%\lim_{n \to \infty}\frac{1}{n} \sum_{k=0}^{n-1} f(\phi_k(z)) = \|f\|_1
%\]
%for almost
% all $z \in \TT$.
% (See, for example, the books of Billingsley~\cite{billingsley} and Halmos~\cite{halmos}.)
%
%
%Thus, for all $w \in H^\infty$ we have that if $\int_{\TT} \log |w| \, dm < 0$ then 
%\[
%w(z)w(\lambda z) \ldots w(\lambda^{n-1}z) \to 0 \qquad \hbox{a.e. on } \TT,
%\]
%since this holds for $w$ outer, and an extra inner factor clearly does not affect the convergence.
%
%{\bf Can we go further?}

The results from Section~\ref{sec:2} can be used to obtain results on the convergence of
iterates of weighted composition operators, although 
there is little published work on the essential spectra of such operators, except in the case of
automorphims $\phi$. Here we
mention  \cite{GZ16, HN15}, and in
 particular, we have the following results from \cite{HN15}, which are clearly
relevant to our analysis, at least in the case that $w$ lies in the disc algebra $A(\DD)$ and is bounded away from $0$ on $\DD$.
\begin{itemize}
\item If $\phi$ is an elliptic automorphism of infinite order with fixed point $\alpha\in \DD$ and $w$ is as above,
then
\[
\sigma(\T)=\sigma_e(\T)= \{\lambda \in \CC: |\lambda|=|w(\alpha)|\}.
\]
\item The same holds if $\phi$ is a parabolic automorphism with fixed point $\alpha \in \TT$.
\item If $\phi$ is a hyperbolic automorphism with attractive/repulsive fixed points $\alpha, \beta$ respectively.
Then 
\[
\sigma(\T)=\sigma_e(\T) = \{\lambda\in \CC: |w(b)|/\phi'(b)^{1/2} \le |\lambda| \le |w(a)|/\phi'(a)^{1/2}\}.
\]
\end{itemize}

Clearly, the essential spectra of weighted composition operators is   a matter for further investigation, but we record here one easy corollary
of the above results.

\begin{cor}
Let $\phi$ be an elliptic automorphism of infinite order with fixed point $\alpha\in \DD$ and $w \in A(\DD)$
bounded away from $0$ on $\DD$. Suppose that $\sup_n \|\T^n\|<\infty$. Then $(\T^n)$
converges uniformly if and only if $|w(\alpha)|<1$.
\end{cor}
\beginpf
This follows directly from Theorem \ref{thm:2.3} on applying the above result from \cite{HN15}.
\endpf

\begin{rem}
In general, suppose that $\phi$ is a finite Blaschke product with $\phi(0)=0$ and that $w \in H^\infty$.
Then for weak (or indeed uniform) convergence we require 
$w(0) \in \{0,1\}$ (also $r_e(\T)<1$ to ensure uniform convergence).

In this situation $C_\phi$ is a contraction, so we have 
the power-boundedness condition $\sup_n \|\T^n\|<\infty$ whenever $\|w\|_\infty \le 1$, and hence also weak convergence.
Indeed, $\|\T^n\| \to 0$ whenever $\|w\|_\infty < 1$.
\end{rem}

\section{Particular classes of operator}
\subsection{Isometries}

If $T$ is an isometry on a Hilbert space, not equal to the identity, then clearly the
sequence $(T^n)$ cannot converge uniformly, or even strongly; however, if $T$ is completely non-unitary then $(T^n)$  will
converge weakly to $0$.

The isometric weighted composition operators on $H^2$ were characterised in \cite{kumar}
as follows.

\begin{prop} $\T$ is an isometry on $H^2$ if and only if $\phi$ is inner, $\|w\|_2=1$ and
$\langle w,w\phi^n\rangle=0$ for all $n \ge 1$.  
\end{prop}

It is worth exploring this further, and   in fact all inner functions $\phi$
can occur as part of an isometry $\T$. 
\begin{prop}
For every inner function $\phi$ there is a weight $w$ such that
$\T$ is an isometry.
\end{prop}
\beginpf
First, if $\phi(0)=0$,  then $C_\phi$ is an isometry, and any inner function $w$ provides an
isometry $\T$.

Next, if $\phi$ is inner and has a zero at $\beta \in \DD$, then we claim
that $w=\theta k_\beta/\|k_\beta\|$ provides an isometry $\T$ for any inner function $\theta$,
where $k_\beta$ is the reproducing kernel at $\beta$.

For $\|w\|_2=1$, and for $n \ge 1$ we have
\[
\langle w,w\phi^n \rangle = \langle \theta k_\beta, \theta k_\beta \phi^n \rangle/\|k_\beta\|^2 = 0,
\]
 since
$\phi(\beta)=0$.

Finally, if $\phi$ is a singular inner function, let 
\[
w= \theta (\phi-\phi(0) )/\|\phi-\phi(0) \|,
\]
where again $\theta$ is inner. We have, for $n \ge 1$,  
\begin{eqnarray*}
&&\langle \theta (\phi-\phi(0) ),\theta (\phi-\phi(0) ) \phi^n \rangle = 
\langle   (\phi-\phi(0) ),   (\phi-\phi(0) ) \phi^n \rangle \\
&=&  \langle 1,  (\phi-\phi(0) )\phi^{n-1}\rangle - \phi(0) \langle 1,  (\phi-\phi(0) )\phi^n \rangle =0,
\end{eqnarray*}
using the fact that $1$ is the reproducing kernel $k_0$. Hence $\T$ is an isometry in this case as well.
\endpf

Thus, unlike in the unweighted case $C_\phi$,
when the isometries require $\phi$ to be inner and $\phi(0)=0$
(see, e.g. \cite{nordgren,CP03}), the Denjoy--Wolff point does not play a role here.
Some further ideas about shifts and unitary operators are given in \cite{matache14}.\\

Isometries on the Bergman space $A^2$ are harder to characterise. Some conditions are given by Zorboska \cite{Z18},
expressed mainly in terms of the measure $\mu$ defined on Borel subsets of $\DD$
by
\beq
\label{eq:defmu}
\mu(A) = \int_{\phi^{-1}(A)} |w(z)|^2 \, dm(z),
\eeq
where now $m$ denotes normalized area measure.
This has the property that
\[
\|\T f\|^2 = \int_\DD |f|^2 \, d\mu(z)
\]
for $f \in A^2$ (see \cite{kumar}).

In the case of isometries, we have the following result, which can be seen as a reformulation of
\cite[Thm. 2.1 (ii)]{Z18}, with a different proof.

\begin{thm}\label{thm:mmu}
The weighted composition operator $\T$ is an isometry on $A^2$ if and only if the
measure $\mu$ defined in \eqref{eq:defmu}
coincides with Lebesgue measure $m$.
\end{thm}
\beginpf
The operator $\T$ is an isometry if and only if
\[
\int_\DD |f(z)|^2 dm(z)= \int_\DD |f(z)|^2 d\mu(z) \qquad (f \in A^2),
\]
so clearly this happens if $\mu=m$. For the converse, 
we introduce a new algebra $\GA \subset C(\overline \DD)$, defined by
\[
\GA= \left\{ \sum_{k=1}^n c_k |p_k|^2: n \in \NN, c_k \in \CC, p_k \in \CC[z] \right \}.
\]
This algebra is self-adjoint, separates points (if $z_1 \ne z_2$, simply choose a 
polynomial that vanishes at $z_1$ but not $z_2$), and contains the constant functions.
Hence, by the Stone--Weierstrass theorem it is dense in $C(\overline \DD)$. Using the
isometry condition  we see that $m$ and $\mu$ induce identical
functionals on $\GA$, hence on $C(\overline\DD)$, and thus $\mu=m$.
\endpf

\begin{rem}
The above methods can be used to characterise isometries on the weighted Bergman
spaces $A^2_\alpha$ for $\alpha>-1$, defined by replacing $m$ by the measure
$(1-|z|^2)^\alpha dm(z)$. We omit the details.
\end{rem}

Zorboska notes that if $\phi$ is a finite Blaschke product of degree $N$, then
$\T$ is isometric on $A^2$ if we take $w(z)=\dfrac{1}{\sqrt N} \phi'(z)$. In fact similar methods
can be used to prove the following more general result.

\begin{prop}
Suppose that $\phi: \DD \to \DD$ is holomorphic and that there is an integer $N \in \NN$ such that
\[
m(\DD \setminus \{z \in \DD: \# \{ \phi^{-1}(\{z\})=N \})=0.
\]
Then $w=c\phi'/\sqrt N$ induces an isometric weighted composition operator $\T$ on $A^2$ for each constant
$c$ with $|c|=1$.
\end{prop}
\beginpf
Let $\Omega$ denote $\{z \in \DD: \# \{ \phi^{-1}(\{z\})=N \}$. Then
\begin{eqnarray*}
\frac{1}{N} \int_\DD |f \circ \phi(z)|^2 |\phi'(z)|^2 \, dm(z) &=& \frac{1}{N} \int_\Omega |f \circ \phi(z)|^2 |\phi'(z)|^2 \, dm(z) \\
&=& \frac{N}{N} \int_\DD |f(s)|^2 dm(s) = \|f\|^2_{A^2}.
\end{eqnarray*}
\endpf

This allows us to construct some interesting new examples of isometries. For example,
let $\psi$ denote the conformal mapping (Riemann mapping) from $\DD$ onto the left semidisc
$\DD_L=\{z \in \DD: \re z < 0\}$ mapping real $z$ to real $z$. Then $\phi:= \psi^4$ is a double cover of
the disc $\DD$ with the exception of the positive real axis, which is covered just once (and $0$ is not
covered at all). Thus $T_{\phi'/\sqrt2, \phi}$ is an isometry of $A^2$.
Explicit formulae for $\psi$ may be found in \cite[p. 64]{CPbook}.\\

However, it is not necessary to use $w=c \phi'/\sqrt{N}$ to obtain Bergman space
isometries. There are many instances when a Blaschke product $\phi$ of degree $N$ has the property that
$\phi \circ \psi = \phi$ for some elliptic automorphism $\psi$ of order $N$ (see, particularly,
\cite[Cor. 2.6]{mortini} for details). In this case we have the following.

\begin{cor}
Let $\phi$ be a Blaschke product of degree $N$ such that $\phi \circ \psi=\phi$,
where $\psi$ is an elliptic automorphism with fixed point $p \in \DD$ and order $N$.
Then the composition operator $\T$ is an isometry on $A^2$ if 
\[
w(z) = c\frac{1}{\sqrt N} \phi'(\psi(z)) \qquad (z \in \DD)
\]
with $|c|=1$.
\end{cor}
\beginpf
Let $A \subset \DD$  be a Borel set and
let $B$ be a Borel set mapped bijectively to $A$ by $\phi$
to within sets of measure $0$
 (this can be obtained by
dividing $\DD$ into $N$ regions, each mapped almost everywhere onto $\DD$). Then
\[
 \int_{\phi^{-1}(A)} |w(z)|^2 \, dm(z) = \int_B \sum_{k=1}^N |w(\psi_k(z))|^2 \, dm(z),
\]
\[
= \frac{1}{N}\int_B \sum_{k=1}^N |\phi'(\psi_{k+1}(z))|^2 \, dm(z)
=  \frac{1}{N}\int_{\phi^{-1}(A)} |\phi'(z)|^2 \, dm(z)
\]
and since this equals $m(A)$ we see that $\T$ is an isometry by applying Theorem \ref{thm:mmu}.
\endpf

\subsection{Hardy spaces of simply-connected domains}

Iterates of (unweighted) composition operators $C_\Phi$
on the right half-plane $\CC_+$ have been studied in \cite{KS18}. We recall a few basic facts about such
operators:
\begin{itemize}
\item For $\Phi: \CC_+ \to \CC_+$ holomorphic, the composition operator $C_\Phi$ on $H^2(\CC_+)$ is bounded if and only if $\Phi(\infty)=\infty$
and its angular derivative satisfies $\Phi'(\infty)>0$. Then the norm, essential norm, and spectral radius of
$C_\Phi$ all equal $\Phi'(\infty)^{1/2}$. (See \cite{EJ12}.)
\item Let $M$ denote the self-inverse conformal mapping $M(z)=(1-z)/(1+z)$ between $\DD$ and $\CC_+$.
Then $C_\Phi$ on $H^2(\CC_+)$ is unitarily equivalent to the operator $\T$ on $H^2(\DD)$,
where $\phi= M \circ \Phi \circ M$ and 
$w(z)= \dfrac{1+\phi(z)}{1+z}$. (See \cite{CP03}.)
\item For bounded $\T$ as above, we have $\DW(\phi)=-1$ and $\phi'(-1)=\Phi'(\infty)$. (See \cite{KS18}.)
\end{itemize}
Some of the conclusions of \cite{KS18} are the following:
\begin{thm}
 For $\Phi: \CC_+ \to \CC_+$ holomorphic, 
 such that $\Phi(\infty)=\infty$ and $0 < \Phi'(\infty) < \infty$, with $\Phi'(\infty) \ne 1$, 
 the following are equivalent:\\
(i)  $C^n_\Phi \to 0$ strongly;\\
(ii) $\DW(\Phi)=\infty$ and $\Phi'(\infty)<1$;\\
(iii) $r_e(C_\Phi)<1$;\\
(iv) $C^n_\Phi \to 0$ uniformly.
\end{thm}

Now let $\Omega$ be a  bounded simply-connected domain with rectifiable (finite-length) boundary,
and $\beta: \DD \to \Omega$ a conformal bijection. Then
let $E^2(\Omega)$ be the Hardy--Smirnoff space of all functions $f: \Omega \to \CC$ such that
the function $Jf: f \mapsto f(\beta(z))\beta'(z)$ lies in $H^2(\DD)$, endowed
with the norm $\|f\|_{E^2(\Omega)}=\|Jf\|_{H^2}$.

The space $E^2(\Omega)$ coincides with Rudin's space $H^2(\Omega)$ of all analytic functions $f$ such
that $|f|^2$  has a harmonic majorant, if and only if there exist constants $a, b > 0$
such that $a \le |\beta'(z)| \le b$ for all $z \in \DD$. We refer to the book of Duren \cite{duren} for
further background on these spaces.

For $\Phi: \Omega \to \Omega$ holomorphic, the composition operator $C_\Phi$ on
$E^2(\Omega)$ is unitarily equivalent to the weighted composition $\T$ on $H^2$,
where $\phi=\beta^{-1}\circ \Phi \circ \beta$, and $w(z)=\left( \beta'(z)/\beta'(\phi(z)) \right)^{1/2}$.
(See, e.g. \cite{kumar}.)

\begin{thm}
Suppose that $\DW(\phi) \in \TT$. Then $(\T^n)$ is unbounded in norm, and thus cannot
converge even weakly.
\end{thm}
\beginpf
Let us write $k_w$ for the reproducing kernel
on $E^2(\Omega)$, so that $\langle f, k_w \rangle=f(w)$ for all $f \in E^2(\Omega)$ and
$w \in \Omega$. Then we have the standard formula  $C^*_\Phi k_w=k_{\Phi(w)}$,
and hence
$(C^*_\Phi)^n k_w= k_{\Phi_n(w)}$ for all $w \in \Omega$.

Fix any $w \in \Omega$, and then the sequence $(\Phi_n(w))=(\beta(\phi_n(\beta^{-1}(w))))$
has at least one accumulation point on the boundary, since $\phi_n(\beta^{-1}(w)) \to \DW(\phi)$.

Now for each point $a \in \partial \Omega$ there are functions in $E^2(\Omega)$ that are unbounded 
near $a$, such as a continuous branch of $(z-a)^{-1/4}$, and it follows that $(k_{\Phi_n(w)})$ is unbounded, and hence $\T$ is not power-bounded.
\endpf

% Some other situations are already covered by the results of Section \ref{sec:3}.

\subsection{Operators on weighted spaces}
\label{sec:6}

In \cite{CP14} the following theorem is given; it is based on the notion of an operator represented by
an infinite lower-triangular matrix. 
For a sequence $d=(d_n)_{n=0}^\infty$ of positive real numbers,
$H^2(d)$ denotes the space of holomorphic functions
$f: z \mapsto \sum_{n=0}^\infty c_n z^n$ such that
\[
\|f\|_{H^2(d)}^2 := \sum_{n=0}^\infty |c_n|^2 d_n^2 < \infty.
\]
Also, for $a \in \DD$, $\psi_a$ is the automorphism defined by \eqref{eq:tnf}.

\begin{thm}
Let $\phi$ be a holomorphic self-mapping of the open unit disc $\DD$,
and let $h \in H^2$ be such that $\T$ is a bounded operator on $H^2$.
Let $(d_n)$ be a decreasing sequence of positive reals such that $H^2(d)$ is
automorphism-invariant. Then $\T$ is also a bounded operator on $H^2(d)$
with 
\[
\|\T\|_{H^2(d)} \le \|C_{\psi_{\phi(0)}}\|_{H^2}\|C_{\psi_{\phi(0)}}\|_{H^2(d)}\|\T\|_{H^2}.
\]
\end{thm}
Easy examples of such   $H^2(d)$ spaces are the weighted Bergman spaces $A^2_\alpha(\DD)$ for $\alpha>-1$ with
$d_n=(n+1)^{-\alpha-1}$ for $n=0,1,2,\ldots$.

As regards convergence of iterates, we may use this theorem to
analyse the case that $\DW(\phi) \in \DD$.

\begin{thm}
Let $\phi$ be a holomorphic self-mapping of the open unit disc $\DD$ with $\DW(\phi)\in \DD$,
and let $h \in H^2$ be such that $\T$ is a bounded operator on $H^2$.
Let $(d_n)$ be a decreasing sequence of positive reals such that $H^2(d)$ is
automorphism-invariant.\\
(i) If $\|\T^n\|_{H^2} \to 0$, then $\|\T^n\|_{H^2(d)} \to 0$.\\
(ii) If $\T^n \to 0$ strongly on $H^2$, then $\T^n \to 0$ strongly on $H^2(d)$.
\end{thm}

\beginpf
(i) We have
\[
\|\T^n\|_{H^2(d)} \le \|C_{\psi_{\phi_n(0)}}\|_{H^2}\|C_{\psi_{\phi_n(0)}}\|_{H^2(d)}\|\T^n\|_{H^2}.
\]
Now $\phi_n(0) \to \DW(\phi)$, and it remains to show that this implies that
the operators $C_{\psi_{\phi_n(0)}}$ are uniformly bounded on $H^2(d)$ -- clearly they are uniformly bounded
on $H^2$. But this follows from \cite[Cor. 2.5]{GP13}, where the bound
\[
\|C_{\psi_a}\|_{H^2(d)} \le \left( \frac{1+|a|}{1-|a|}\right)^{\Lambda/2}
\]
is given, with $\Lambda$ being a constant depending on $d$.\\

(ii) Given that the iterates $\T^n$ converge strongly to $0$ on $H^2$ it follows from the fact that $(d_n)$ is bounded that $\|\T^n f \|_{H^2(d)} \to 0$ for all polynomials $f$, which form a dense subset of $H^2(d)$.
But the sequence of operators $(\T^n)$ is uniformly bounded
on
$H^2$,  and hence also on $H^2(d)$ by arguments similar to those in (i),
and so $\|\T^n f\|_{H^2(d)} \to 0$ for all $f \in H^2(d)$.
\endpf

\end{document}